\documentclass[12pt]{article}
\usepackage{graphicx}
\usepackage{amsmath,amssymb,amsthm,amsfonts}
\newtheorem{thm}{Theorem}[section]

\newtheorem{rmk}[thm]{Remark}

\newcommand{\R}{{\mathbb{R}}}

\newcommand{\1}{\partial}
\newcommand{\2}{\overline}
\newcommand{\3}{\varepsilon}

\def\ni{\noindent}

\begin{document}
\title{Another proof for the removable singularities\\ 
of the heat equation}
\author{Kin Ming Hui\\
Institute of Mathematics, Academia Sinica,\\
Nankang, Taipei, 11529, Taiwan, R. O. C.}
\date{Sept 2, 2009}
\smallbreak \maketitle
\begin{abstract}
We give two different simple proofs for the removable singularities
of the heat equation in $(\Omega\setminus\{x_0\})\times (0,T)$ where
$x_0\in\Omega\subset\R^n$ is a bounded domain with $n\ge 3$. We also give
a necessary and sufficient condition for removable singularities
of the heat equation in $(\Omega\setminus\{x_0\})\times (0,T)$ 
for the case $n=2$. 
\end{abstract}

\vskip 0.2truein

Key words: removable singularities, heat equation

Mathematics Subject Classification: Primary 35B65 Secondary 35K55, 35K05, 
35K20
\vskip 0.2truein
\setcounter{equation}{0}
\setcounter{section}{0}

%\section{}
\setcounter{equation}{0}
\setcounter{thm}{0}

Singularities of solutions of partial differential equations appear in 
many problems. For example singularities appears in the study of the
solutions of the harmonic map \cite{SU} and the harmonic map heat flow 
\cite{GGT}. In \cite{SY} S.~Sato and E.~Yanagida studied the solutions
for a semilinear parabolic equation with moving singularities. 
Singularities of solutions also appears in the study of hyperbolic 
partial differential equations \cite{S} and in the study of the touchdown 
behavior of the micro-electromechanical systems equation \cite{GG1},
\cite{GG2}, \cite{GG3}. 

It is interesting to find the necessary and sufficient condition
for the solutions of the equations to have removable singularities. 
In \cite{Hs2} S.Y.~Hsu proved the following theorem.

\begin{thm}
Let $n\ge 3$ and let $0\in\Omega\subset\R^n$ be a domain. Suppose $u$ is a 
solution of the heat equation 
\begin{equation}
u_t=\Delta u
\end{equation}
in $(\Omega\setminus\{0\})\times (0,T)$. 
Then $u$ has removable singularities at $\{0\}\times (0,T)$ if and only if
for any $0<t_1<t_2<T$ and $\delta\in (0,1)$ there exists $\2{B_{R_0}(0)}
\subset\Omega$ depending on $t_1$, $t_2$ and $\delta$, such that 
\begin{equation}
|u(x,t)|\le\delta |x|^{2-n}
\end{equation}
for any $0<|x|\le R_0$ and $t_1\le t\le t_2$. 
\end{thm}
The proof in \cite{Hs2} is based on the Green function estimates of 
\cite{H1} and a careful analysis of the behavior of the 
solution near the singularities using Dehamel principle. In this paper 
we will use the Schauder estimates for heat equation \cite{F}, \cite{LSU}, 
and the technique of \cite{DK} and \cite{Hs1} to give two different 
simple proofs of the above result. We also obtain the following result 
for the solution of the heat equation in $2$-dimension.

\begin{thm}
Let $0\in\Omega\subset\R^2$ be a domain. Suppose $u$ is a 
solution of the heat equation in $(\Omega\setminus\{0\})\times (0,T)$. 
Then $u$ has removable singularities at $\{0\}\times (0,T)$ if and only if
for any $0<t_1<t_2<T$ and $\delta\in (0,1)$ there exists $\2{B_{R_0}(0)}
\subset\Omega$ depending on $t_1$, $t_2$ and $\delta$, such that 
\begin{equation}
|u(x,t)|\le\delta (\log (1/|x|))^{-1}
\end{equation}
for any $0<|x|\le R_0$ and $t_1\le t\le t_2$. 
\end{thm}
\begin{rmk}
Note that the function $\log |x|$ satisfies the heat equation in 
$(\R^2\setminus\{0\})\times (0,\infty)$ but it has non-removable 
singularities on $\{0\}\times (0,\infty)$ and it does not satisfy (3). 
Hence (3) is sharp.
\end{rmk}

We start with some definitions. For any set $A$ we let $\chi_A$ be the 
characteristic function of the set $A$. Let $0\in\Omega\subset\R^n$ 
be a bounded domain. We say that a solution $u$ of the heat equation (1)
in $(\Omega\setminus\{0\})\times (0,T)$ has removable singularities 
at $\{0\}\times (0,T)$ if there exists a classical solution $v$ of (1) in 
$\Omega\times (0,T)$ such that $u=v$ in $(\Omega\setminus\{0\})\times 
(0,T)$. For any $R>0$ let $B_R=B_R(0)=\{x:|x|<R\}\subset\R^n$.

\vskip 0.1truein
{\ni{\it Proof of Theorem 1}:} Suppose $u$ has removable singularities 
at $\{0\}\times (0,T)$. By the same argument as the proof in section 3
of \cite{Hs2}
for any $0<t_1<t_2<T$ and $\delta\in (0,1)$ there exists $\2{B}_{R_0}
\subset\Omega$ depending on $t_1$, $t_2$ and $\delta$, such that (2) 
holds.

Suppose (2) holds. Then for any $0<t_1<t_2<T$ and $\delta\in (0,1)$ 
there exists $\2{B}_{R_0}\subset\Omega$ depending on $t_1$, $t_2$ and 
$\delta$, such that (2) holds for any $0<|x|\le R_0$ and $t_1\le t\le t_2$. 

For any $0<|x|\le R_0$, let
\begin{equation}
w(y,s)=u(|x|y,|x|^2s)\quad\forall 0<|y|\le R_0/|x|,t_1/|x|^2\le s
\le t_2/|x|^2. 
\end{equation}
Then $w$ is a solution of (1) in $(\2{B}_1\setminus\{0\})\times
(|x|^{-2}t_1,|x|^{-2}t_2)$. By (2),
\begin{equation}
|w(y,s)|\le\delta (|x||y|)^{2-n}\quad\forall 0<|y|\le R_0/|x|,t_1/|x|^2\le s
\le t_2/|x|^2. 
\end{equation}
Let $t_1<t_3<t_2$. Then
\begin{equation}
\frac{t_3}{|x|^2}-\frac{t_1}{|x|^2}\ge\frac{t_3-t_1}{R_0^2}
\end{equation}
By the parabolic Schauder estimates \cite{F}, \cite{LSU}, (5) and (6),
there exists a constant $C_1>0$ such that
\begin{equation}
|\nabla w(y,s)|\le C_1\sup_{\tiny\begin{array}{c}
1/2\le |z|\le 1\\
|x|^{-2}t_1\le\tau\le |x|^{-2}t_2\end{array}}
w(z,\tau)\le C_2\delta|x|^{2-n}
\end{equation}
holds for any $2/3\le |y|\le 3/4$, $t_3/|x|^2\le s\le t_2/|x|^2$ where 
$C_2=2^{n-2}C_1$. By (4) and (7),
\begin{align}
&|\nabla u(z,t)|\le C_2\delta|x|^{1-n}\quad\forall |z|=\frac{3}{4}|x|, 
0<|x|\le R_0, t_3\le t\le t_2\nonumber\\
\Rightarrow\quad&|\nabla u(z,t)|\le C_2\delta |z|^{1-n}
\quad\forall |z|\le\frac{3}{4}R_0, t_3\le t\le t_2.
\end{align}
Let $R_1=3/(4R_0)$.
We will now use a modification of the proof of Lemma 2.3 of \cite{DK}
and Lemma 2.1 of \cite{Hs1} to complete the argument. We will first 
show that $u$ satisfies (1) in $\Omega\times (t_1,t_2)$ in the 
distribution sense. Since $u$ satisfies (1) in $(\Omega\setminus\{0\})
\times (0,T)$, for any $0<\3<R_1$ and $\eta\in C_0^{\infty}
(\Omega\times (0,T))$ we have
\begin{align}
\int_{\Omega\setminus B_{\3}}u\eta\,dx\biggr |_{t_3}^{t_2}
=&\int_{t_3}^{t_2}\int_{\Omega\setminus B_{\3}}u\eta_t\,dxdt
-\int_{t_3}^{t_2}\int_{\Omega\setminus B_{\3}}\nabla u
\cdot\nabla\eta\,dxdt\nonumber\\
&\qquad -\int_{t_3}^{t_2}\int_{\1 B_{\3}}\eta \frac{\1 u}{\1 n}\,
d\sigma dt
\end{align} 
where $\1 u/\1 n$ is the derivative of $u$ with respect to the unit
outward normal at $\1 B_{\3}$. By (8),
\begin{equation*}
\limsup_{\3\to 0}\left|\int_{t_3}^{t_2}\int_{\1 B_{\3}}\eta 
\frac{\1 u}{\1 n}\,d\sigma dt\right|\le C_2\delta(t_2-t_3)|\1 B_1|
\|\eta\|_{L^{\infty}}
\end{equation*}
Since $\delta>0$ is artibrary, there holds
\begin{equation}
\lim_{\3\to 0}\int_{t_3}^{t_2}\int_{\1 B_{\3}}\eta \frac{\1 u}{\1 n}\,
d\sigma dt\,dxdt=0.
\end{equation}
By (8) and the Lebesgue dominated convergence theorem,
\begin{equation}
\lim_{\3\to 0}\int_{t_3}^{t_2}\int_{\Omega\setminus B_{\3}}\nabla u
\cdot\nabla\eta\,dxdt
=\int_{t_3}^{t_2}\int_{\Omega}\nabla u\cdot\nabla\eta\,dxdt
\end{equation}
Letting $\3\to 0$ in (9), by (10) and (11) there holds
\begin{equation}
\int_{\Omega}u\eta\,dx\biggr |_{t_3}^{t_2}
=\int_{t_3}^{t_2}\int_{\Omega}u\eta_t\,dxdt
-\int_{t_3}^{t_2}\int_{\Omega}\nabla u
\cdot\nabla\eta\,dxdt\quad\forall t_3\in (t_1,t_2).
\end{equation}
Hence $u$ is a distribution solution of (1) in $\Omega\times (t_1,t_2)$.
By (2) for any $1\le p<\frac{n}{n-2}$ there exists a constant $C_p'>0$
such that 
\begin{equation}
\sup_{t_1\le t\le t_2}\int_{B_{R_0}}u(x,t)^p\,dx\le C_p'
\end{equation}
By (12) and (13) and an argument similar to the proof of \cite{KZ} and
section 1 of \cite{H2} $u\in L_{loc}^{\infty}(B_{R_0}\times (t_1,t_2))$. 
We now let $v$ be the solution of 
\begin{equation}
\left\{\begin{aligned}
&v_t=\Delta v\qquad\qquad\qquad\text{ in }
B_{R_1}\times (t_3,t_2)\\
&\frac{\1 v}{\1 n}(x,t)=\frac{\1 u}{\1 n}(x,t)\quad\mbox{ on }
\1 B_{R_1}\times (t_3,t_2)\\
&v(x,t_3)=u(x,t_3)\qquad\mbox{ in }B_{R_1}.
\end{aligned}\right.
\end{equation}
For any $0\le h\in C_0^{\infty}(B_{R_1})$ and $t_3<t\le t_2$ let $\eta$ be the 
solution of
\begin{equation}
\left\{\begin{aligned}
&\eta_t+\Delta\eta=0\quad\text{ in }
B_{R_1}\times (t_3,t)\\
&\frac{\1\eta}{\1 n}(x,t)=0\quad\mbox{ on }\1 B_{R_1}\times (t_3,t)\\
&\eta (x,t)=h(x)\quad\mbox{in }B_{R_1}.
\end{aligned}\right.
\end{equation}
By the maximum principle, 
\begin{equation}
0\le\eta\le \|h\|_{L^{\infty}}\quad\mbox{ in }B_{R_1}\times (t_3,t).
\end{equation}
Then by (14) and (15),
\begin{equation}
\begin{aligned}
\left.\int_{B_{R_1}\setminus B_{\3}}(u-v)\eta\,dx\right|_{t_3}^t
=&\int_{t_3}^t\int_{B_{R_1}\setminus B_{\3}}[(u-v)\eta_t+(u-v)_t\eta]
\,dxdt\\
=&\int_{t_3}^t\int_{B_{R_1}\setminus B_{\3}}[(u-v)\eta_t+\Delta 
(u-v)\eta]\,dxdt\\
=&\int_{t_3}^t\int_{B_{R_1}\setminus B_{\3}}(u-v)(\eta_t+\Delta\eta)\,dxdt
-\int_{t_3}^t\int_{\1 B_{\3}}\eta\frac{\1}{\1 n}(u-v)\,d\sigma dt\\
&\qquad+\int_{t_3}^t\int_{\1 B_{\3}}(u-v)\frac{\1\eta}{\1 n}\,d\sigma dt\\
=&-\int_{t_3}^t\int_{\1 B_{\3}}\eta\frac{\1}{\1 n}(u-v)\,d\sigma dt
+\int_{t_3}^t\int_{\1 B_{\3}}(u-v)\frac{\1\eta}{\1 n}\,d\sigma dt.
\end{aligned}
\end{equation}
By (2),
\begin{equation}
\left|\int_{t_3}^t\int_{\1 B_{\3}}(u-v)\frac{\1\eta}{\1 n}\,d\sigma dt
\right|\le C\3\to 0\quad\mbox{ as }\3\to 0.
\end{equation}
By (8) and (16),
\begin{equation}
\limsup_{\3\to 0}\left|\int_{t_3}^t\int_{\1 B_{\3}}\eta
\frac{\1}{\1 n}(u-v)\,d\sigma dt\right|\le C\delta.
\end{equation}
Since $\delta>0$ is arbitrary, by (19) there holds
\begin{equation}
\lim_{\3\to 0}\left|\int_{t_3}^t\int_{\1 B_{\3}}\eta
\frac{\1}{\1 n}(u-v)\,d\sigma dt\right|=0.
\end{equation}
Letting $\3\to 0$ in (17), by (18) and (20), 
\begin{equation}
\int_{B_{R_1}}(u-v)(x,t)h(x)\,dx
=\int_{B_{R_1}}(u-v)(x,t_3)\eta (x,t_3)\,dx=0.
\end{equation}
We now choose a sequence of functions $h_i\in C_0^{\infty}(B_{R_1})$
converging to $\chi_{\{u>v\}}$ a.e. $x\in B_{R_1}$ as $i\to\infty$. 
Putting $h=h_i$ in (21) and letting $i\to 0$,
\begin{equation}
\int_{B_{R_1}}(u-v)_+(x,t)\,dx=0\quad\forall t_3<t\le t_2.
\end{equation}
By interchanging the role of $u$ and $v$ we get
\begin{equation}
\int_{B_{R_1}}(v-u)_+(x,t)\,dx=0\quad\forall t_3<t\le t_2.
\end{equation}
Hence by (22) and (23),
\begin{align}
&\int_{B_{R_1}}|v-u|(x,t)\,dx=0\quad\forall t_3<t\le t_2\nonumber\\
\Rightarrow\quad&u(x,t)=v(x,t)\quad\forall 0<|x|\le R_1, t_3<t\le t_2.
\end{align}
Hence $u$ has removable singularities on $\{0\}\times (t_3,t_2)$.
Since $0<t_1<t_3<t_2<T$ is arbitrary, $u$ has removable singularities 
on $\{0\}\times (0,T)$ and the theorem follows.

\vskip 0.1truein
{\ni{\it Proof of Theorem 2}:} Theorem 2 follows by an argument very
similar to the proof of Theorem 1 but with (3) replacing (2) in the 
argument.

\vskip 0.1truein
{\ni{\it An alternate proof of Theorems 1 and 2}:} We will show that 
when (2) (respectively (3)) holds, then $u$ has removable singularities
at $\{0\}\times (0,T)$. Suppose (2) holds if $n\ge 3$ and (3) holds
if $n=2$. We first observe that by the previous argument for any 
$0<t_1<t_2<T$ $u$ satisfies (12) and $u\in L_{loc}^{\infty}(\Omega
\times (0,T))$. Let $\2{B}_{R_1}\subset\Omega$ and let $w$ be the 
solution of 
\begin{equation*}
\left\{\begin{aligned}
&w_t=\Delta w\quad\mbox{ in }B_{R_1}\times (t_1,t_2)\\
&w=u\qquad\mbox{ on }\2{B}_{R_1}\times\{t_1\}\cup
\1 B_{R_1}\times (t_1,t_2).
\end{aligned}\right.
\end{equation*}
By the maximum principle, 
\begin{equation}
\|w\|_{L^{\infty}}\le\|u\|_{L^{\infty}(B_{R_1}\times (t_1,t_2))}<\infty.
\end{equation}
For any $\3>0$, let 
\begin{equation*}
w_{\3}=\left\{\begin{aligned}
&w-u+\3|x|^{2-n}\qquad\quad\,\,\,\mbox{if }n\ge 3\\
&w-u+\3\log (R_1/|x|)\quad\mbox{ if }n=2.
\end{aligned}\right.
\end{equation*}
Then $w_{\3}$ satisfies
\begin{equation*}
\left\{\begin{aligned}
&w_{\3,t}=\Delta w_{\3}\quad\mbox{ in }(B_{R_1}\setminus\{0\})
\times (t_1,t_2)\\
&w_{\3}\ge u\qquad\quad\mbox{ on }\1 B_{R_1}\times (t_1,t_2)\cup
\2{B}_{R_1}\times\{t_1\}.
\end{aligned}\right.
\end{equation*}
By (2), (3), and (25) there exists a constant $0<r_0<R_1$ such that
$$
w_{\3}\ge 0\quad\mbox{ on }\1 B_{r_1}\times [t_1,t_2] 
$$
for all $0<r_1\le r_0$. By the maximum principle in $(B_{R_1}\setminus
B_{r_1})\times (t_1,t_2)$,
\begin{align}
&w_{\3}\ge 0\quad\mbox{ in }(B_{R_1}\setminus B_{r_1})\times (t_1,t_2)
\nonumber\\
\Rightarrow\quad&\left\{\begin{aligned}
&w-u+\3|x|^{2-n}\ge 0\qquad\quad\forall r_1\le |x|\le R_1, 
t_1\le t\le t_2\quad\,\mbox{if }n\ge 3\\
&w-u+\3\log (R_0/|x|)\ge 0\quad\forall r_1\le |x|\le R_1, 
t_1\le t\le t_2\quad\mbox{ if }n=2
\end{aligned}\right.\nonumber\\
\Rightarrow\quad&w\ge u\qquad\quad\forall 0<|x|\le R_1, t_1\le t\le t_2
\quad\mbox{ as }r_1\to 0,\3\to 0.
\end{align}
Similarly by considering the function
\begin{equation*}
v_{\3}=\left\{\begin{aligned}
&w-u-\3|x|^{2-n}\qquad\quad\,\mbox{if }n\ge 3\\
&w-u-\3\log (R_1/|x|)\quad\mbox{ if }n=2
\end{aligned}\right.
\end{equation*}
and applying the maximum principle and letting $\3\to 0$ we get
\begin{equation}
w\le u\quad\forall 0<|x|\le R_1, t_1\le t\le t_2.
\end{equation}
By (26) and (27) we get (24) and Theorem 2 and Theorem 3 follows.

\end{document}